\documentclass[12pt]{amsart}
\usepackage{graphicx} 
\usepackage{amsthm}
\usepackage{amsmath}
\usepackage{amssymb}

\newtheorem{definition}{Definition}

\usepackage{xcolor}
\usepackage[colorlinks]{hyperref}
\usepackage[colorinlistoftodos]{todonotes}
\usepackage{verbatim}
\usepackage{comment}
\usepackage{subcaption}
\usepackage{listings}

\def\Z{\mathbb{Z}}

\definecolor{Nat}{rgb}{0.61, 0.77, 0.89}

\title{A note on girth-diameter cages}

\author[Araujo-Pardo]{Gabriela Araujo-Pardo}
\address{Gabriela Araujo-Pardo, Instituto de Matem\'aticas-Campus Juriquilla, Universidad Nacional Aut\'onoma de M\'exico, M\'exico.}
\email{garaujo@im.unam.mx}

\author[Conder]{Marston Conder}
\address{Marston Conder, Department of Mathematics, University of Auckland, New Zealand}
\email{m.conder@auckland.ac.nz}

\author[Garc\'ia-Colin]{Natalia Garc\'ia-Col\'in}
\address{Natalia Garc\'ia-Col\'in, D\'epartement d'Informatique, Universit\'e Libre de Bruxelles and Department of Statistical Learning, ScaDS.AI Leipzig}
\email{natalia.garcia.colin@ulb.be}

\author[Kiss]{Gy\"orgy Kiss}
\address{Gy\"orgy Kiss, Department of Geometry and HUN-REN-ELTE Geometric and Algebraic Combinatorics,
Research Group, Eötvös Loránd University, Hungary, and
Faculty of Mathematics, Natural Sciences and Information Technologies,
University of Primorska, Slovenia}
\email{gyorgy.kiss@ttk.elte.hu}

\author[Leemans]{Dimitri Leemans}
\address{Dimitri Leemans, D\'epartement de Math\'ematique, Universit\'e libre de Bruxelles, and Department of Statistical Learning, ScaDS.AI Leipzig}
\email{leemans.dimitri@ulb.be}

\date{\today}

\begin{document}
\maketitle

\begin{abstract}
In this paper we introduce a problem closely related to the Cage Problem and the Degree Diameter Problem.
For integers $k\geq 2$, $g\geq 3$ and $d\geq 1$, we define a $(k;\, g,d)$-graph to be a $k$-regular graph with girth $g$ and diameter $d$. We denote by $n_0(k;\,g,d)$ the smallest possible order of such a graph, and, if such a graph exists, we call it a $(k;g,d)$-cage.
In particular, we focus on $(k;\,5,4)$-graphs. We show that $n_0(k;\,5,4) \geq k^2+k+2$ for all $k$, and report on the determination of all $(k;\,5,4)$-cages for $k=3, 4$ and $5$ and examples with $k = 6$, and describe some examples of $(k;\,5,4)$-graphs which prove that $n_0(k;\,5,4) \leq 2k^2$ for infinitely many values of $k$.
\end{abstract}

\section{Introduction}

The {\em order\/} of a graph is the number of its vertices, its {\em diameter\/} is the maximum distance between a pair of its vertices, and its {\em girth\/} is the length of its smallest circuit.  A graph is $k$-{\em regular\/} if each of its vertices has exactly $k$ neighbours. 

The well-known {\em Cage Problem\/} involves finding the $k$-regular graphs of girth $g$ with smallest possible order $n(k,g)$, for a given pair  $(k,g)$ of integers with $k\geq 2$ and $g\geq 3$. A regular graph with these properties is called a {\em $(k,g)$-cage,} or simply a {\em cage}. Cages were introduced by Tutte \cite{T47} in 1947, and the Cage Problem has been widely studied from the time when Erd\"os and Sachs  \cite{ES63} proved their existence in 1963. A complete survey about this topic and its relevance can be found in \cite{EJ13}. 

The equally well known {\em Degree-Diameter Problem\/} involves finding the $k$-regular graphs of diameter $d$ with largest possible order, for a given pair  $(k,d)$ of integers with $k\geq 2$ and $d\geq 1$. 
In this case what is known as the `Moore bound' states that the largest order is at most $1+k + k(k-1)+\cdots + k(k-1)^{d-1}$. It is easy to see that graphs attaining this upper bound are cages with odd girth $g = 2d+1$. Conversely, if $g$ is odd then this is also a lower bound for the order of a $k$-regular graph with girth $g$, when $d=(g-1)/2$.  For a complete review of the Degree-Diameter problem, we refer the reader to~\cite{MS2013}. 

Motivated by the above background, in this paper we describe a variation of the Cage Problem by considering a lower bound for the order of regular graphs with given girth and given diameter, which we call {\em girth-diameter cages\/} (or simply $gd$-{\em cages\/}).   

Specifically, for given integers $k,g$ and $d$ with $k\geq 2$, $g\geq 3$ and $d\geq 1$, we define a $(k;g,d)$-graph to be a $k$-regular graph with girth $g$ and diameter $d$, and we denote by $n_0(k;\,g,d)$ the smallest possible order of such a graph (if one exists).  A necessary condition for existence is $\lfloor \frac{g}{2}\rfloor\leq d$ (since $g \le 2d+1$), but our main interest in such graphs is in cases where $\lfloor \frac{g}{2}\rfloor \leq d\leq g$. 

In this note, we study the first pair of parameters that we considered interesting, namely those with girth $g = 5$ and diameter $d = 4$, that is, we study $(k;\,5,4)$-graphs.

For any $(k;\,5,4)$-graph $G$ of order $n$, let $r$ and $c$ be two vertices at distance $4$ from each other.  Then the neighbourhoods $N(r)$ and $N(c)$ of $r$ and $c$ must each consist of $k$ vertices, with $N(c)\cap N(r)= \emptyset$, and the remaining $n-2k-2$ vertices of $G$ form a set $M$ that includes all neighbours of vertices in $N(r)$ apart from $r$ and all neighbours of vertices in $N(c)$ apart from $c$.  Hence $|M| \ge k(k-1)$, and it follows that $|V(G)| \ge 2+2k+k(k-1) = k^2 + k + 2$.  We may call this number $k^2 + k + 2$ the Moore bound for $n_0(k;\,5,4)$, and denote it by $M(k;\,5,4)$.

Note that it is easy to obtain a generic bound for any given integers $k,g$ and $d$ with $k\geq 2$, $g\geq 5$ and $d\geq 4$. 

\section{$(k;\,5, 4)$-cages for $k=3,4,5$ and $6$.} \label{sec:gd-cages}

For degree $k=3,$ we found an example of a $(3;\,5,4)$-graph attaining the Moore bound $n = M(3;\,5,4) = 3^2 + 3 + 2 = 14$, and by an easy computational search we found that up to isomorphism there is just one other. These two $(3;\,5,4)$-cages are depicted in Figure~\ref{Banffrican}, with the $19$ black edges in common.  One of them contains also the two red edges (and its automorphism group has order $12$), while the other one contains the two green edges (and its automorphism group has order $4$). They can be found in House of Graphs~\cite{hog} with identification numbers 1000 and 50487.  

\begin{center}
${}$\\[-18pt]
 \begin{figure}[htbp]
\includegraphics[scale=0.8]{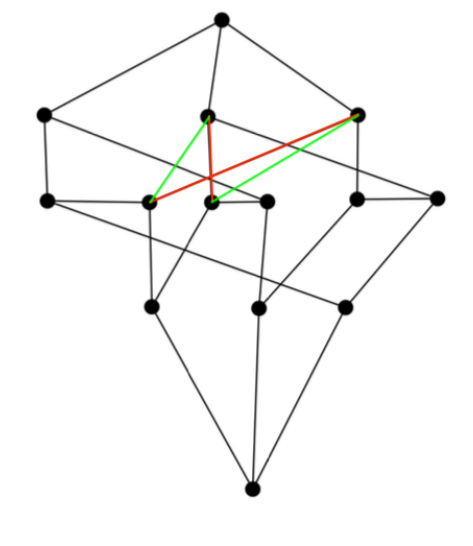}
${}$\\[-18pt]
\caption{The two graphs with degree 3, diameter 4 and girth 5}\label{Banffrican}
\end{figure}   
\end{center}

Similarly for degree $k=4,$ we found that there are exactly four non-isomorphic graphs attaining the Moore bound $M(4;\,5,4) = 4^2 + 4 + 2 = 22$, with automorphism groups of orders $1$, $2$, $4$ and $8$.
These graphs are available in House of Graphs with identification numbers 50459, 49991, 49992 and 49993 respectively.

For degree $k=5,$ we found using a more advanced computational search that there are exactly seven non-isomorphic graphs that attain the Moore bound $M(5;\,5,4) = 5^2 + 5 + 2 = 32$, with 
automorphism groups of orders $4$, $4$, $10$, $10$, $48$, $64$ and $1920$.  
Two of the latter seven $(5;\,5,4)$-cages appeared already (in a different context) in~\cite{GD2022}, namely those with automorphism groups of orders $48$ and $1920$. 

A similar computational search has so far also produced two non-isomorphic $6$-valent graphs that attain the Moore bound $M(6;\,5,4) = 6^2 + 6 + 2 = 44$, with automorphism groups of orders $40$ and $240$.


\section{$(k;\, 5,4)$-graphs from the Levi graphs of 
biaffine planes.}


Some small $(k;\, 5,4)$-graphs graphs can be obtained using an amalgamation method on the Levi graphs (incidence graphs) of biaffine planes, 
as we explain below.

\begin{definition}
Let $\Pi _q$ be a finite projective plane of order $q$. A biaffine plane is obtained from 
$\Pi _q$ by choosing a point-line pair $(P,\ell )$, deleting $P$, deleting $\ell$, all the lines incident with $P$ and all the points belonging to $\ell$. If the point-line pair is incident in $\Pi _q$, then we say that the biaffine plane is of type $1$, and otherwise we say it is of type $2$.
\end{definition}

The Levi graph of $\Pi _q$ is a $(q+1,6)$-cage and its diameter is $3.$ The Levi graphs 
of biaffine planes are $(q;\, 6,4)$-graphs of orders $2q^2$ and $2(q^2-1),$
respectively. Two vertices are at distance $4$ if and only if both of them correspond either to points on a deleted line or to lines in the same parallel class. 

The first time that this type of construction was used to find graphs of girth $5$ (as far as the authors are aware), is in a paper \cite{Br-5} by Brown, who proved that if $q\geq 5$ is a prime power, then $n(q+2,5)\leq 2q^2.$ 

More sophisticated but similar construction methods have been presented by Funk \cite{funk}, Abreu \emph{et al} \cite{AABL11}, and Abajo and her co-authors in a series of papers \cite{AABB17, ABBM19, AB21}. In order to construct these graphs, the authors of these papers applied the so-called `amalgamation' technique.
By inserting some new edges between vertices in the incidence graph of the biaffine plane that correspond either to two points of a deleted line or to two lines of a parallel class, they obtained regular graphs with degree $q+a$ for certain values of $a\geq 2$, and with diameter $4$ when the diameter of the amalgamated subgraph is at least $4$. For example, this happens with a
$q$-cycle in the incidence graph of a biaffine plane of order $q$, for $q\geq 7$.

The exact orders of these $(k;\, 5,4)$-graphs are a little smaller than $2k^2$, but generally larger than $M(k;\, 5, 4)=k^2+k+2$.  

As mentioned before, Araujo-Pardo and Leemans \cite{GD2022} used a biaffine plane of type $1$ to construct a $(q+1,5,4)$-graph of order $2q^2$ for $q=4$, and this particular graph attains the lower bound $M(5;\, 5,4) = 32$.  Then recently Araujo-Pardo, Kiss and 
Porups\'anszki \cite{AKP23} generalised the latter graph for any $2$-power $q=2^r,$ and presented a simple geometric construction for $(q+2;\, 5,4)$-graphs whose order is $2q^2-2$, using a biaffine plane of type $2$.

\section{The `middle' graph}

For construction of the various  examples we found for degree $k= 5$ and $6$, we strongly utilised properties of the subgraph $H$ induced by the `middle' vertex set $M = V(G) \setminus (\{r\} \cup \{c\} \cup N(r) \cup N(c))$ described earlier, when $G$ is a graph attaining the Moore bound $M(k;\,5,4)$. In particular:

\begin{enumerate}
    \item[{\rm (a)}] $H$ has order $k(k-1)$;
    \item[{\rm (b)}] $H$ is regular with degree $k-2$; 
    \item[{\rm (c)}] $H$ has girth at least $5$, so contains no circuits of length $3$ or $4$; 
    \item[{\rm (d)}] $H$ has diameter at least $3$; and 
    \item[{\rm (e)}] the vertex-set of $H$ can be labelled with ordered pairs $(i, j)$ with $i \ne j \in \Z_k$, such that 
 if $(i,j)$ is adjacent to $(i', j')$ then $i\neq i'$ and $j\neq j'$,
and two vertices $(i,j)$ and $(i',j')$ with $i = i'$ or $j = j'$ must lie at distance at least 3 from each other.
\end{enumerate} 

The first three properties above are obvious. To verify the last two, let us introduce some extra notation. 
Denote the $k$ vertices in $N(r)$ as $r_1, \ldots r_k$, and those in $N(c)$ as $c_1, \ldots c_k$, 
and for $1\le i \le k$ define $R_i = N(r_i) \setminus \{r\}$ and $C_i = N(c_i) \setminus \{c\}.$ 

Because $H$ contains no $3$-cycles, we find that $R_i \cap R_j = \emptyset$ and  $C_i \cap C_j = \emptyset$ whenever $i \neq j$, and because $H$ contains no $4$-cycles, also $|C_i \cap R_i| \leq 1$ for all $i,j$. Hence we can think of each $R_i$ as a row and each $C_j$ as a column of the $k \times k$ integer grid corresponding to $\Z_k \times Z_k$ with its diagonal $\{(j,j) : j \in \Z_k\}$ removed, and relabel the vertices of $H$ as the pairs $(i, j)$ with $ i, j \in \Z_k$ such that $i \ne j$. 

Next, there must be at least one pair of vertices in $H$ at distance at least $3$ from each other, for otherwise a path of length $2$ between vertices $(i,j)$ and $(i', j)$ in $H$ would form a $4$-cycle in $G$ when taken together with the edges from $c_j$ to $(i,j)$ and $(i', j)$. This line of argument also implies that two vertices that share the first or the second coordinate must lie at distance at least $3$ from each other.
\medskip

The examples we found (and described earlier in Section \ref{sec:gd-cages}) have the following properties: 

\begin{enumerate} 
\item[$\bullet$]  When $k = 3$ the subgraph $H$ of order $6$ consists of three non-incident edges;
\item[$\bullet$]  When $k = 4$ the subgraph $H$ is a $2$-regular graph of order $12$, and indeed in all four examples, it is a single cycle, with girth $12$ and diameter $6$; 
\item[$\bullet$]  When $k = 5$ the subgraph $H$  is a $3$-regular graph of order $20$, with diameter $4$ and girth $5$ in two cases,
diameter $5$ and girth $5$ in three cases, and diameter $4$ and girth $6$ in the other two cases.  
\end{enumerate} 

\medskip

We now present the following:

\smallskip\noindent
{\bf Proposition.}
A graph $H$ with the properties  (a) to (e) listed above can be extended to a $(k;\,5, 4)$-cage. 

\begin{proof}
Suppose that $H$ is a graph satisfying the conditions (a) to (e). 
Then the set of vertices and edges of a new graph $G$ can be constructed by adding the $2+2k$ vertices $r, r_1, \ldots r_k, c, c_1, \ldots , c_k$ and the edges $\{r, r_i\}$, $\{r_i, (i,j)\}$, $\{c, c_j\}$, $\{c_j, (i,j)\}$ for $i, j \in \Z_k$, and we immediately find that $G$ is a $k$-regular graph with $k(k+1)+2$ vertices. All that remains is to check that the diameter of $G$ is $4$ and its girth is $5$.

First it is clear that every $r_i$ and every $c_j$ lies at distance at most $3$ from each of $r$ and $c$, and every vertex $(i,j)$ lies at distance $2$ from each of $r$ and $c$, and hence at distance at most $3$ from every $r_{i'}$ and every $c_{j'}$, and at distance at most $4$ from every other middle vertex $(i',j')$.  It now follows easily that every two vertices lie at distance at most $4$ from each other, so the diameter is $4$.

Finally, suppose that some circuit of length $3$ or $4$ occurs in $G$.  Then since the girth of $H$ is at least $5$, at least one of the vertices of that circuit lies in $H$ and at least one lies in $G - H$.  The latter vertex cannot be $r$ or $c$, so that circuit must contain some $r_i$ or some $c_j$, but neither of $r$ and $c$. Hence the circuit contains two vertices $(i,j)$ and $(i',j')$ with $i = i'$ or $j = j'$, but then by property (e) those two vertices lie at distance at least $3$ from each other in $H$, so the circuit must contain another vertex of $G - H$ either adjacent to both $(i,j)$ and $(i,j')$, or adjacent to both $(i,j)$ and $(i',j)$, which is impossible. Hence the girth of $G$ is $5$.
\end{proof}

\section{Future work}
Using the above Proposition, the search for $(k;\,5,4)$-cages may be reduced to the search for suitable `middle' graphs satisfying the conditions (a) to (e) given in the previous section.  Some of us are continuing to search for possible constructions of such graphs that would work for other values of the degree $k\geq 6$.

Finally, we remark that the middle graph $H$ is a $(k-2,g)$-graph where $g\geq 5$ and the graph is of order
$$k(k-1)=(k-2)^2+3(k-2)+2=M(k-2;\, 5,4)+2(k-2),$$
and that for $k\geq 14$ the order of the record holder $(k-2,5)$-graphs is greater than this number (see \cite{EJ13} and Table 1 in  \cite{ABBM19}). Hence it would seem to be quite a challenge to find $(k;\,5,4)$-cages of order $k^2+k+2$ for large $k$.

\section{Acknowledgements}
This research was initiated at the BIRS workshop 23w5125 `Extremal Graphs arising from Designs and Configurations'. The authors are very grateful to the organisers and to BIRS. 
Gabriela Araujo-Pardo acknowledges  support by PAPIIT-M\'exico under grant IN101821. 
Marston Conder acknowledges support from New Zealand's Marsden Fund (project UOA2030).
Dimitri Leemans acknowledges support from an Action de Recherche Concert\'ee grant from the Communaut\'e Française -- Wallonie Bruxelles and the Fonds National de la Recherche Scientifique de Belgique. 
Gy\"orgy Kiss acknowledges support from the Hungarian National Research, Development and Innovation Office OTKA grant no. SNN 132625. 

\bibliographystyle{plain}

\end{document}